\documentclass[12pt]{article}

\usepackage{amsmath}
\usepackage{amsthm}
\usepackage{amssymb}
\usepackage{latexsym}
\usepackage{graphicx, pst-plot, pst-node, pst-text, pst-tree}

\newtheorem{theorem}{Theorem}[section]
\newtheorem{proposition}[theorem]{Proposition}

\newtheorem{corollary}[theorem]{Corollary}

\newtheorem{conjecture}[theorem]{Conjecture}

\newcounter{todocounter}

 \makeatletter
\newfont{\footsc}{cmcsc10 at 8truept}
\newfont{\footbf}{cmbx10 at 8truept}
\newfont{\footrm}{cmr10 at 10truept}
\title{Grid classes and the Fibonacci dichotomy for restricted permutations}
\author{Sophie Huczynska\thanks{Supported by a Royal Society Dorothy Hodgkin Research Fellowship.}\ \ and Vincent Vatter\thanks{Supported by EPSRC grant GR/S53503/01.}\\[-5pt]
\small School of Mathematics and Statistics\\[-5pt]
\small University of St.\ Andrews\\[-5pt]
\small St.\ Andrews, Fife, Scotland\\[-5pt]
\small \texttt{\{sophieh, vince\}@mcs.st-and.ac.uk}\\[-5pt]
\small \texttt{http://turnbull.mcs.st-and.ac.uk/\~{}\{sophieh,
vince\}}}

\date{\today \\[6pt]
    \begin{flushleft}
    \end{flushleft}
}

\begin{document}
\maketitle

\newcommand{\st}{\operatorname{st}}
\newcommand{\Av}{\operatorname{Av}}
\newcommand{\C}{\mathcal{C}}
\newcommand{\vecv}{{\bf v}}
\newcommand{\x}{{\bf x}}
\newcommand{\y}{{\bf y}}
\newcommand{\zpm}{0/\mathord{\pm} 1}
\newcommand{\Grid}{\operatorname{Grid}}
\renewcommand{\sb}{\mathcal{SB}}

\begin{abstract}
We introduce and characterise grid classes, which are natural generalisations of other well-studied permutation classes.  This characterisation allows us to give a new, short proof of the Fibonacci dichotomy: the number of permutations of length $n$ in a permutation class is either at least as large as the $n$th Fibonacci number or is eventually polynomial.
\end{abstract}

\section{Introduction}

A permutation $\pi$ of $[n]$
\footnote{Here
$[n]=\{1,2,\ldots,n\}$ and, more generally, for $a,b \in
\mathbb{N}$ ($a<b$), the interval $\{a,a+1,\ldots,b\}$ is
denoted by $[a,b]$, the interval $\{a+1,a+2,\ldots,b\}$ is denoted by
$(a,b]$, and so on.}
contains the permutation $\sigma$ of $[k]$ ($\sigma\le\pi$) if $\pi$ has a subsequence of length $k$
in the same relative order as $\sigma $.
For example, $\pi=391867452$ (written in list, or one-line notation)
contains $\sigma=51342$, as can be seen by considering the subsequence
$91672$ ($=\pi(2),\pi(3),\pi(5),\pi(6),\pi(9)$).
A {\it permutation class\/} is a downset of permutations under this order, or in other words, if $\C$ is a permutation class, $\pi\in\C$, and $\sigma\le\pi$, then $\sigma\in\C$.
We shall denote by $\C_n$ ($n \in \mathbb{N}$)
the set $\C \cap S_n$, i.e.\ those permutations in $C$ of length
$n$.  
Recall that
an {\it antichain} is a set of pairwise incomparable elements.
For any permutation class $\C$, there is a unique (and possibly infinite) antichain $B$ such that $\C=\Av(B)=\{\pi: \beta \not \leq\pi\mbox{ for all } \beta \in B\}$. This antichain $B$ is called the {\it basis} of $\C$.
Permutation classes arise naturally in a variety of disparate fields, ranging from the analysis of sorting machines (dating back to Knuth~\cite{knuth1}, who proved that a permutation is stack-sortable if and only if it lies in the class $\Av(231)$) to the study of Schubert varieties (see, e.g., Lakshmibai and Sandhya~\cite{ls:smooth}).

The Stanley-Wilf Conjecture, recently proved by Markus and Tardos~\cite{mt:swc}, states that all permutation classes except the set of all permutations have at most exponential growth, i.e., for every class $\C$ with a nonempty basis, there is a constant $K$ so that $|\C_n|<K^n$ for all $n$.  Less is known regarding the exact enumeration of permutation classes.  Natural enumerative questions include:
\begin{enumerate}
\item[(i)] Which permutation classes are finite?
\item[(ii)] Which permutation classes are enumerated by a polynomial?
\item[(iii)] Which permutation classes have rational generating functions?  (We refer to $\sum |\C_n|x^n$ as the generating function of $\C$.)
\item[(iv)] Which permutation classes have algebraic generating functions?
\item[(v)] Which permutation classes have $P$-recursive enumeration?
\end{enumerate}
The answer to the first question on this list follows easily from the Erd\H os-Szekeres Theorem\footnote{{\bf The Erd\H os-Szekeres Theorem~\cite{es:acpig}.} 
Every permutation of length $n$ contains a monotone subsequence of length at least $\sqrt{n}$.}: the class $\Av(B)$ is finite if and only if $B$ contains both an increasing permutation and a decreasing permutation.
The answer to the second question is provided in this paper.  Questions (iii)--(v) remain unanswered.

\bigskip\noindent{\bf Downsets of vectors.} Perhaps the simplest interesting context in which to study downsets is finite vectors of nonnegative integers, and in this context there is also a polynomial enumeration result which we shall employ in our proofs.  Let
$\x=(x_1,x_2,\dots,x_m),\y=(y_1,y_2,\dots,y_m)\in\mathbb{N}^m$ for
some $m$.  We write $\x\le\y$ if $x_i\le y_i$ for all $i\in[m]$.
This order is often called the {\it product order\/}.  The weight
of the vector $\x$, denoted $\|\x\|$, is the sum of the entries of
$\x$.

\begin{figure}
\begin{center}
\psset{xunit=0.01in, yunit=0.01in} \psset{linewidth=1\psxunit}
\begin{pspicture}(0,0)(150,150)
\psaxes[dy=10, Dy=1, dx=10, Dx=1, tickstyle=bottom,
showorigin=false, labels=none](0,0)(140,140)
\multips(0,0)(10,0){15}{\pscircle*{4.0\psxunit}}
\multips(0,10)(0,10){14}{\pscircle*{4.0\psxunit}}
\multips(10,10)(10,0){14}{\pscircle*{4.0\psxunit}}
\multips(10,20)(0,10){13}{\pscircle*{4.0\psxunit}}
\multips(20,20)(0,10){13}{\pscircle*{4.0\psxunit}}
\multips(30,20)(0,10){13}{\pscircle*{4.0\psxunit}}
\pscircle*(40,20){4.0\psxunit} \pscircle*(50,20){4.0\psxunit}
\pscircle*(60,20){4.0\psxunit}
\multips(40,30)(0,10){5}{\pscircle*{4.0\psxunit}}
\pscircle(70,20){4.0\psxunit} \pscircle(50,30){4.0\psxunit}
\pscircle(40,80){4.0\psxunit}
\end{pspicture}
\end{center}
\caption{The plot of downset in $\mathbb{N}^2$; the elements
of the class are drawn with solid circles, while the elements of
the basis are drawn with hollow circles.}
\end{figure}
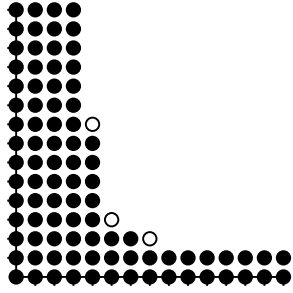

\begin{theorem}\label{vector-poly}
Let $\C$ denote a downset in $\mathbb{N}^m$.  For
sufficiently large $n$, the number of vectors in $\C$ of weight
$n$ is given by a polynomial.
\end{theorem}

Stanley~\cite{s:E2546} posed Theorem~\ref{vector-poly} as a
{\it Monthly\/} problem in 1976 and offered two solutions.  One of these solutions is elementary while the other follows from viewing the number of vectors in question as a Hilbert function.

\bigskip\noindent{\bf Downsets of other objects.} Downsets of other combinatorial objects have been extensively studied, and other polynomial enumeration results are known.  These have often been established by ideas analogous to the grid classes of matchings we use.

For example, downsets of graphs w.r.t.\ the induced subgraph ordering that are closed under isomorphism are called {\it hereditary properties\/}.  Let $\mathcal{P}$ denote a hereditary property, and let $\mathcal{P}_n$ denote the set of graphs in $\mathcal{P}$ with vertex set $[n]$.  Scheinerman and Zito~\cite{sz:hered} proved that $|\mathcal{P}_n|$ either has polynomial growth (meaning that $|\mathcal{P}_n|=\Theta(n^k)$ for some $k$) or $|\mathcal{P}_n|$ has at least exponential growth.  Balogh, Bollob{\'a}s, and Weinreich~\cite{bbw:growth2000} later showed that polynomial growth hereditary properties are enumerated exactly by a polynomial for large $n$.  Their proof of this result uses ``canonical properties,'' which are quite like our grid classes of matchings.

Moving to a more general context, Pouzet and Thi\'ery~\cite{pt:poly} study polynomial growth (although not exact polynomial enumeration) for certain downsets of relational structures.  While summarising their work would take us too far afield, we remark first that permutations can be viewed as relational structures\footnote{E.g., $\pi\in S_n$ can be taken to correspond to the relational structure on $[n]$ with two linear orders, $<$ and $\prec$, where $<$ is the normal ordering of $[n]$ and $i\prec j\iff \pi(i)<\pi(j)$.} and second that the grid classes of matchings we use essentially correspond to their concept of ``monomorphic decompositions into finitely many parts.''

\section{Grid classes}\label{grid-classes}

\subsection{The skew-merged permutations}\label{skew-merged}

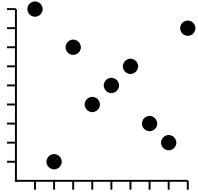
\begin{figure}[t]
\begin{center}
\psset{xunit=0.01in, yunit=0.01in} \psset{linewidth=1\psxunit}
\begin{pspicture}(0,0)(100,100)
%\psline[linecolor=gray](55,0)(55,100)
%\psline[linecolor=gray](0,55)(100,55)
\psaxes[dy=10, Dy=1, dx=10, Dx=1, tickstyle=bottom,
showorigin=false, labels=none](0,0)(90,90)
\pscircle*(10,90){4.0\psxunit} \pscircle*(20,10){4.0\psxunit}
\pscircle*(30,70){4.0\psxunit} \pscircle*(40,40){4.0\psxunit}
\pscircle*(50,50){4.0\psxunit} \pscircle*(60,60){4.0\psxunit}
\pscircle*(70,30){4.0\psxunit} \pscircle*(80,20){4.0\psxunit}
\pscircle*(90,80){4.0\psxunit}
\end{pspicture}
\end{center}
\caption{The plot of the skew-merged permutation
$917456328$.}\label{fig-skew-merged}
\end{figure}

We begin with an example of a grid class.  A permutation is said to be {\it skew-merged\/} if it is the union of an increasing subsequence and a decreasing subsequence.  For example, the permutation shown in Figure~\ref{fig-skew-merged} is skew-merged.  Stankova~\cite{stankova:fs} was the first to find the basis of this class.  Later, K\'ezdy, Snevily, and Wang~\cite{ksw:incdec} observed that the basis follows easily from F\"oldes and Hammer's characterisation of split graphs%
\footnote{A graph $G$ is {\it split\/} if its vertices can be
partitioned into a disjoint union $V(G)=V_1\uplus V_2$ s.t.\ $G[V_1]$
is complete and $G[V_2]$ is edgeless.  F\"oldes and
Hammer proved that a graph is split if and only if it does not
contain $K_2\uplus K_2$, $C_4$, or $C_5$ as induced subgraphs.} in
\cite{fh:sg}.

\begin{theorem}[Stankova~\cite{stankova:fs}; K\'ezdy, Snevily, and Wang~\cite{ksw:incdec}; and Atkinson~\cite{a:skewmerged}]\label{skewmerged}
The skew-merged permutations are $\Av(2143, 3412)$.
\end{theorem}

Atkinson~\cite{a:skewmerged} gave another proof of
Theorem~\ref{skew-merged} and solved the enumeration problem for
the skew-merged permutations; he showed that their generating
function is given by
$$
\frac{1-3x}{(1-2x)\sqrt{1-4x}}.
$$

K\'ezdy, Snevily, and Wang~\cite{ksw:incdec} studied one
generalization of skew-merged permutations, the class of
permutations which can be partitioned into $r$ increasing
subsequences and $s$ decreasing subsequences.  Grid classes
provide a different generalization.

\subsection{Definitions}

First an important warning: when discussing grid classes,
we index matrices beginning from the lower left-hand corner, and we
reverse the rows and columns; for example $M_{3,2}$ denotes for us
the entry of $M$ in the $3$rd column from the left and $2$nd row
from the bottom.  Below we include a matrix with its entries
labeled:
$$
\left(
\begin{footnotesize}
\begin{array}{rrr}
(1,2)&(2,2)&(3,2)\\
(1,1)&(2,1)&(3,1)
\end{array}
\end{footnotesize}
\right).
$$

\bigskip

\begin{figure}[t]
\begin{center}
\psset{xunit=0.01in, yunit=0.01in} \psset{linewidth=1\psxunit}
\begin{pspicture}(0,0)(100,100)
\psline[linecolor=gray](55,0)(55,100)
\psline[linecolor=gray](0,55)(100,55) \psaxes[dy=10, Dy=1, dx=10,
Dx=1, tickstyle=bottom, showorigin=false, labels=none](0,0)(90,90)
\pscircle*(10,90){4.0\psxunit} \pscircle*(20,10){4.0\psxunit}
\pscircle*(30,70){4.0\psxunit} \pscircle*(40,40){4.0\psxunit}
\pscircle*(50,50){4.0\psxunit} \pscircle*(60,60){4.0\psxunit}
\pscircle*(70,30){4.0\psxunit} \pscircle*(80,20){4.0\psxunit}
\pscircle*(90,80){4.0\psxunit}
\end{pspicture}
\end{center}
\caption{A gridding of the skew-merged permutation
$917456328$.}\label{fig-skew-merged-gridded}
\end{figure}
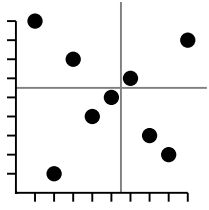

Roughly, the {\it grid class\/} of a matrix $M$ is the set of all
permutations that can be divided in a prescribed manner (dictated
by $M$) into a finite number of blocks, each containing a monotone subsequence.  We have already
introduced the best-studied grid class, the skew-merged
permutations.  We previously defined them as the permutations that
can be written as the union of an increasing subsequence and a
decreasing subsequence.  As a grid class, the skew-merged
permutations can be defined as the permutations that can be
divided into four monotonic blocks, two increasing and two
decreasing, as indicated in Figure~\ref{fig-skew-merged-gridded},
and our notation for this class is
$$
\Grid\left(\begin{footnotesize}\begin{array}{rr}-1&1\\1&-1\end{array}\end{footnotesize}\right),
$$
but before reaching that point we need to introduce some notation.

Given a permutation $\pi\in S_n$ and sets $A,B\subseteq[n]$,
we write
$\pi(A\times B)$ for the subsequence of $\pi$ with indices
from $A$ which has values in $B$.  For example, applying this
operation to the permutation shown in
Figure~\ref{fig-skew-merged-gridded}, we get
$$
917456328([5]\times[5])=1,4,5,
$$
and this (increasing) subsequence gives the points in the
lower left-hand box of Figure~\ref{fig-skew-merged-gridded}.  The
increasing subsequence in the upper right-hand box is
$$
917456328([6,9]\times[6,9])=6,8,
$$
while the decreasing subsequence in the lower right-hand box is
$$
917456328([6,9]\times[5])=3,2.
$$

Now suppose that $M$ is a $t\times u$ matrix (meaning, in the
notation of this paper, that it has $t$ columns and $u$ rows).
An {\it $M$-gridding\/} of the permutation $\pi\in S_n$ is a pair
of sequences $1=c_1\le c_2\le\cdots\le c_{t+1}=n+1$ (the column
divisions) and $1=r_1\le r_2\le\cdots\le r_{u+1}=n+1$ (the row
divisions) such that for all $k\in[t]$ and $\ell\in[u]$,
$\pi([c_k,c_{k+1})\times[r_\ell,r_{\ell+1}))$ is:
\begin{itemize}
\item increasing if $M_{k,\ell}=1$, \item decreasing if
$M_{k,\ell}=-1$, \item empty if $M_{k,\ell}=0$.
\end{itemize}
We define the {\it grid class of $M$\/}, written $\Grid(M)$, to be
the set of all permutations that possess an $M$-gridding.  We say
that $\pi$ is {\it $t\times u$-griddable\/} if it is $M$-griddable
for some $t\times u$ matrix $M$.

A class $\C$ is said to be {\it $t\times u$-griddable\/} if every
permutation in $\C$ is $t\times u$-griddable, and it is said to be
{\it griddable\/} if it is $t\times u$-griddable for some
$t,u\in\mathbb{N}$.  Note that all griddable classes lie in some
particular grid class (suppose that $\C$ is $t\times u$ griddable
and take a larger matrix $M$ containing every $t\times u$ matrix,
then $\C$ lies in $\Grid(M)$).

Two special types of grid classes have been extensively studied.
One type is the profile classes of Atkinson~\cite{a:rp}, which in
our language are grid classes of permutation matrices.  Another
example of grid classes are the $W$-classes introduced by
Atkinson, Murphy, and Ru\v{s}kuc~\cite{amr:pwocsop}, which are the
grid classes of $\zpm$ row vectors.

Atkinson, Murphy, and Ru\v{s}kuc~\cite{amr:pwocsop} introduced $W$-classes in their study of partially well-ordered (pwo)\footnote{Recall that a partially ordered set is said to be {\it partially well-ordered (pwo)\/} if it contains neither an infinite properly decreasing sequence
nor an infinite antichain.} permutation classes, and proved that
grid classes of $\zpm$ row vectors are pwo.  This result
does not extend to arbitrary grid classes, i.e., some grid classes
contain infinite antichains, e.g., there is an infinite antichain of
skew-merged permutations.
In order to characterise the pwo grid
classes, we associate a graph to each grid class.  For any
$t\times u$ matrix $M$ we construct the bipartite graph $G(M)$
with vertices $x_1$, $x_2$, $\dots$, $x_t$ and $y_1$, $y_2$,
$\dots$, $y_u$ and edges $x_ky_\ell$ precisely when
$M_{k,\ell}\neq 0$.  For example, the bipartite graph of a vector
is a star together with isolated vertices, while the bipartite
graph of
$\left(\begin{footnotesize}\begin{array}{rr}-1&1\\1&-1\end{array}\end{footnotesize}\right)$
is a cycle with $4$ vertices.  The pwo properties of a grid class depend only on its graph.

\begin{theorem}[Murphy and Vatter~\cite{profile}]\label{prof-pwo}
The grid class of $M$ is pwo if and only if $G(M)$ is a forest.
\end{theorem}

\subsection{The characterisation of griddable classes}

It appears surprisingly difficult to compute the basis of
$\Grid(M)$ when $M$ is neither a vector nor a permutation matrix.
Waton [private communication] has computed the bases of
$\Grid(M)$ for all $2\times 2$ matrices $M$, but we know of
no such results for larger matrices.
In particular, the following remains a conjecture.

\begin{conjecture}
All grid classes are finitely based.
\end{conjecture}

We instead take a coarser approach here and ask only for a
characterisation of the griddable classes, that is, the
permutation classes that lie in some grid class.

It will prove useful to have an alternative interpretation of
griddability.  We say that the permutation $\pi\in S_n$ can be
{\it covered by $s$ monotonic rectangles\/} if there are
$[w_1,x_1]\times[y_1,z_1]$,$\dots$,$[w_s,x_s]\times[y_s,z_s]\subseteq[n]\times[n]$
such that
\begin{itemize}
\item for each $i\in[s]$, $\pi([w_i,x_i]\times[y_i,z_i])$ is
monotone, and \item $\displaystyle\bigcup_{i\in[s]}
[w_i,x_i]\times[y_i,z_i]=[n]\times[n]$.
\end{itemize}
Note that we allow these rectangles to intersect.  By definition
every $t\times u$-griddable permutation can be covered by $tu$
monotonic rectangles.  The following proposition gives the
other direction.

\begin{proposition}\label{grid-coverings}
Every permutation that may be covered by $s$ monotonic rectangles
is $(2s-1)\times (2s-1)$-griddable.
\end{proposition}
\begin{proof}
Suppose that $\pi\in S_n$ is covered by the $s$ monotonic
rectangles $[w_1,x_1]\times[y_1,z_1]$, $\dots$,
$[w_s,x_s]\times[y_s,z_s]\subseteq[n]\times[n]$.  Define the
indices $c_1,c_2,\dots,c_{2s}$ and $r_1,r_2,\dots,r_{2s}$ by
\begin{eqnarray*}
\{c_1\le c_2\le\cdots\le c_{2s}\}&=&\{w_1,x_1,w_2,x_2,\dots,w_s,x_s\},\\
\{r_1\le r_2\le\cdots\le
r_{2s}\}&=&\{y_1,z_1,y_2,z_2,\dots,y_s,z_s\}.
\end{eqnarray*}
Since these rectangles cover $\pi$, we must have $c_1=r_1=1$ and
$c_{2s}=r_{2s}=n$.  Now we claim that these sets of indices form
an $M$-gridding of $\pi$ for some $2s-1\times 2s-1$ matrix $M$.

To prove this claim it suffices to show that
$\pi([c_k,c_{k+1}]\times[r_\ell,r_{\ell+1}])$ is monotone for
every $k,\ell\in[2s-1]$, since we can then construct the matrix
$M$ based on whether this subsequence is increasing or decreasing.
Because the rectangles given cover $\pi$, the point $(c_k,r_\ell)$
lies in at least one rectangle, say $[w_m,x_m]\times[y_m,z_m]$.
Thus $c_k\ge w_m$ and $r_\ell\ge y_m$ and, because of the ordering
of the $c$'s and $r$'s, we have $c_{k+1}\le x_m$ and
$r_{\ell+1}\le z_m$.  Therefore
$[c_k,c_{k+1}]\times[r_\ell,r_{\ell+1}]$ is contained in
$[w_m,x_m]\times[y_m,z_m]$ and so
$\pi([c_k,c_{k+1}]\times[r_\ell,r_{\ell+1}])$ is monotone.
\end{proof}

With this new interpretation of griddability established, we need only two more definitions before characterising the griddable classes.  Given two permutations $\pi\in S_m$ and $\sigma\in S_n$, we define their {\it direct sum\/}, written $\pi\oplus\sigma$ by
$$
(\pi\oplus\sigma)(i)
=
\left\{
\begin{array}{ll}
\pi(i)&\mbox{if $i\in [m]$,}\\
\sigma(i-m)+m&\mbox{if $i\in[m+n]\setminus[m]$,}
\end{array}
\right.
$$
and similarly define their {\it skew sum\/} by
$$
(\pi\ominus\sigma)(i)
=
\left\{
\begin{array}{ll}
\pi(i)+n&\mbox{if $i\in [m]$,}\\
\sigma(i-m)&\mbox{if $i\in[m+n]\setminus[m]$.}
\end{array}
\right.
$$

\begin{theorem}\label{griddable-characterisation}
A permutation class is griddable if and only if it does not
contain arbitrarily long direct sums of $21$ or skew sums of $12$.
\end{theorem}
\begin{proof}
If a permutation class does contain arbitrarily long direct sums
of $21$ or skew sums of $12$, then it is clearly not griddable.

For the other direction, let $\C$ be a permutation class that does
not contain $\ominus^{a+1} 12$ or $\oplus^{b+1} 21$.  We show
by induction on $a+b$ that there is a function $f(a,b)$ so that
every permutation in $\C$ can be covered by $f(a,b)$ monotonic
rectangles, and thus we will be done by
Proposition~\ref{grid-coverings}.

First note that if either $a$ or $b$ is $0$ then $\C$ can only
contain monotone permutations, so we can set $f(a,0)=f(0,b)=1$.
The next case is $a+b=2$, and since we may assume that $a,b\neq
0$, we have $a=b=1$.  Thus $\C$ contains neither $\ominus^2
12=3412$ nor $\oplus^2 21=2143$, so $\C$ is a subclass of the
skew-merged permutations and thus every permutation in $\C$ may be
covered by $4$ monotonic rectangles and we may take $f(1,1)=4$.

%\definecolor{gray}{rgb}{.75, .75, .75}

\begin{figure}
\begin{center}
\begin{tabular}{ccccccc}

\psset{xunit=0.02in, yunit=0.02in} \psset{linewidth=0.5\psxunit}
\begin{pspicture}(0,0)(50,50)
\psframe[linecolor=white,fillstyle=vlines,hatchangle=45,hatchcolor=gray](0,0)(20,20)
\psline[linecolor=gray,linestyle=dashed,linewidth=1\psxunit](20,40)(20,0)
\psline[linecolor=gray,linestyle=dashed,linewidth=1\psxunit](0,20)(40,20)
\psaxes[dy=10, Dy=1, dx=10, Dx=1, tickstyle=bottom,
showorigin=false, labels=none](0,0)(50,50)
\pscircle*(10,30){2\psxunit} \pscircle*(20,40){2\psxunit}
\pscircle*(30,10){2\psxunit} \pscircle*(40,20){2\psxunit}
\end{pspicture}

&\rule{10pt}{0pt}&

\psset{xunit=0.02in, yunit=0.02in} \psset{linewidth=0.5\psxunit}
\begin{pspicture}(0,0)(50,50)
\psframe[linecolor=white,fillstyle=vlines,hatchangle=45,hatchcolor=gray](20,0)(50,30)
\psline[linecolor=gray,linestyle=dashed,linewidth=1\psxunit](10,30)(50,30)
\psline[linecolor=gray,linestyle=dashed,linewidth=1\psxunit](20,40)(20,0)
\psaxes[dy=10, Dy=1, dx=10, Dx=1, tickstyle=bottom,
showorigin=false, labels=none](0,0)(50,50)
\pscircle*(10,30){2\psxunit} \pscircle*(20,40){2\psxunit}
\pscircle*(30,10){2\psxunit} \pscircle*(40,20){2\psxunit}
\end{pspicture}

&\rule{10pt}{0pt}&

\psset{xunit=0.02in, yunit=0.02in} \psset{linewidth=0.5\psxunit}
\begin{pspicture}(0,0)(50,50)
\psframe[linecolor=white,fillstyle=vlines,hatchangle=45,hatchcolor=gray](0,20)(30,50)
\psline[linecolor=gray,linestyle=dashed,linewidth=1\psxunit](30,10)(30,50)
\psline[linecolor=gray,linestyle=dashed,linewidth=1\psxunit](0,20)(40,20)
\psaxes[dy=10, Dy=1, dx=10, Dx=1, tickstyle=bottom,
showorigin=false, labels=none](0,0)(50,50)
\pscircle*(10,30){2\psxunit} \pscircle*(20,40){2\psxunit}
\pscircle*(30,10){2\psxunit} \pscircle*(40,20){2\psxunit}
\end{pspicture}

&\rule{10pt}{0pt}&

\psset{xunit=0.02in, yunit=0.02in} \psset{linewidth=0.5\psxunit}
\begin{pspicture}(0,0)(50,50)
\psframe[linecolor=white,fillstyle=vlines,hatchangle=45,hatchcolor=gray](30,30)(50,50)
\psline[linecolor=gray,linestyle=dashed,linewidth=1\psxunit](10,30)(50,30)
\psline[linecolor=gray,linestyle=dashed,linewidth=1\psxunit](30,10)(30,50)
\psaxes[dy=10, Dy=1, dx=10, Dx=1, tickstyle=bottom,
showorigin=false, labels=none](0,0)(50,50)
\pscircle*(10,30){2\psxunit} \pscircle*(20,40){2\psxunit}
\pscircle*(30,10){2\psxunit} \pscircle*(40,20){2\psxunit}
\end{pspicture}

\\[4pt]

(i)&&(ii)&&(iii)&&(iv)

\end{tabular}
\end{center}
\caption{The regions of $\pi$ referred to in the proof of
Theorem~\ref{griddable-characterisation}.}\label{characterisation-fig}
\end{figure}
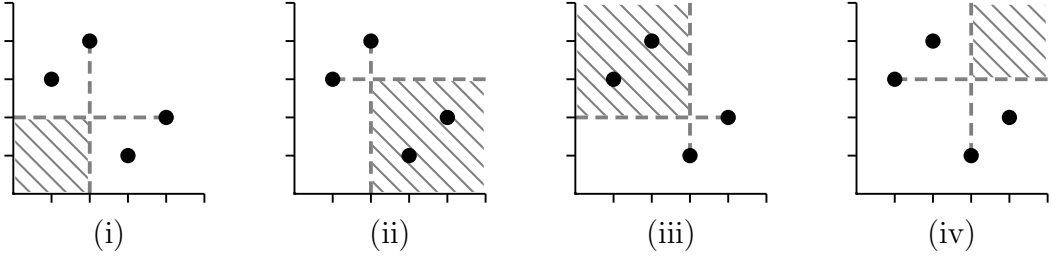

By symmetry and the cases we have already handled, we may assume
that $a\ge 2$ and $b\ge 1$.  Let $\pi\in\C_n$ be a
$3412$-containing permutation (if there are no such permutations,
then we are done by induction) and suppose that
$\pi(i_1)\pi(i_2)\pi(i_3)\pi(i_4)$ is in the same relative order
as $3412$ where $1\le i_1<i_2<i_3<i_4\le n$.  By induction we have
the following (see Figure~\ref{characterisation-fig} for an
illustration of these regions):
\begin{enumerate}
\item[(i)] $\pi([i_2]\times[\pi(i_4)])$ avoids $\ominus^{a+1} 12$
and $\oplus^{b} 21$ so it can be covered by $f(a,b-1)$ monotonic
rectangles, \item[(ii)] $\pi([i_2,n]\times[\pi(i_1)])$ avoids
$\ominus^{a} 12$ and $\oplus^{b+1} 21$ so it can be covered by
$f(a-1,b)$ monotonic rectangles, \item[(iii)]
$\pi([i_3]\times[\pi(i_4),n])$ avoids $\ominus^{a} 12$ and
$\oplus^{b+1} 21$ so it can be covered by $f(a-1,b)$ monotonic
rectangles, and \item[(iv)] $\pi([i_3,n]\times[\pi(i_1),n])$
avoids $\ominus^{a+1} 12$ and $\oplus^{b} 21$ so it can be covered
by $f(a,b-1)$ monotonic rectangles.
\end{enumerate}
Because the four regions in (i)--(iv) cover $\pi$, it may be
covered by $2f(a-1,b)+2f(a,b-1)$ monotonic rectangles.
Furthermore, the $3412$-avoiding permutations in $\C$ may be
covered by $f(1,b)\le f(a-1,b)$ monotonic rectangles by induction,
so we may take $f(a,b)=2f(a-1,b)+2f(a,b-1)$, completing the proof.
\end{proof}

\subsection{The enumeration of grid classes}

To date only scattered results are known about the enumeration of
grid classes and their subclasses.  The only general results are the following two.

\begin{theorem}[Atkinson~\cite{a:rp}]\label{perm-mat-poly}
If $M$ is a permutation matrix, then $\Grid(M)$ and all its
subclasses have eventually polynomial enumeration.
\end{theorem}

\begin{theorem}[Albert, Atkinson, and Ru\v{s}kuc~\cite{aar:regular}]\label{star-rational}
If $G(M)$ is a star, then $\Grid(M)$ and all its subclasses have
rational (and readily computable) generating functions.
\end{theorem}

It is very tempting to speculate that the enumerative properties
of a grid class depend only on its graph\footnote{For example:
\begin{conjecture}\label{forest-rational}
If $G(M)$ is a forest then $\Grid(M)$ and all its subclasses have
rational generating functions.
\end{conjecture}}.  Our contribution to
this suspicion is to show (in Theorem~\ref{grid-poly}) that when
$G(M)$ is a matching\footnote{We take a {\it matching\/} to be a
graph without incident edges, i.e., a graph with maximum
degree $1$.} then $\Grid(M)$ and all its subclasses have
eventually polynomial enumeration, thus generalising
Theorem~\ref{perm-mat-poly}.  For brevity, we refer to such
classes as the {\it grid classes of matchings\/}.

\begin{theorem}\label{grid-poly}
If the permutation class $\C$ lies in the grid class of a matching
then there is a polynomial $p(n)$ so that $|\C_n|=p(n)$ for all
sufficiently large $n$.
\end{theorem}
\begin{proof}
Let $M$ be a $t\times u$ matrix whose graph is a matching, let
$\C$ be a subclass of $\Grid(M)$, and let $\pi\in\C$.  We define
the {\it greedy $M$-gridding\/} of $\pi$ to be the gridding given
by $1=c_1\le c_2\le\cdots\le c_{t+1}=n+1$ (the column divisions)
and $1=r_1\le r_2\le\cdots\le r_{u+1}=n+1$ (the row divisions)
where for each $k$, $c_k$ is chosen so as to maximise
$c_1+c_2+\cdots+c_k$.  Because $G(M)$ is a matching, this uniquely
defines the $r$'s.

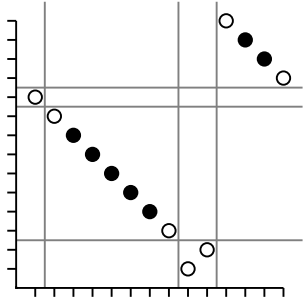
\begin{figure}
\begin{center}
\psset{xunit=0.01in, yunit=0.01in} \psset{linewidth=1\psxunit}
\begin{pspicture}(0,0)(150,150)
\psline[linecolor=gray](15,0)(15,150)
\psline[linecolor=gray](85,0)(85,150)
\psline[linecolor=gray](105,0)(105,150)
\psline[linecolor=gray](0,25)(150,25)
\psline[linecolor=gray](0,95)(150,95)
\psline[linecolor=gray](0,105)(150,105) \psaxes[dy=10, Dy=1,
dx=10, Dx=1, tickstyle=bottom, showorigin=false,
labels=none](0,0)(140,140) \pscircle(10,100){4.0\psxunit}
\pscircle(20,90){4.0\psxunit} \pscircle*(30,80){4.0\psxunit}
\pscircle*(40,70){4.0\psxunit} \pscircle*(50,60){4.0\psxunit}
\pscircle*(60,50){4.0\psxunit} \pscircle*(70,40){4.0\psxunit}
\pscircle(80,30){4.0\psxunit} \pscircle(90,10){4.0\psxunit}
\pscircle(100,20){4.0\psxunit} \pscircle(110,140){4.0\psxunit}
\pscircle*(120,130){4.0\psxunit} \pscircle*(130,120){4.0\psxunit}
\pscircle(140,110){4.0\psxunit}
\end{pspicture}
\end{center}
\caption{A greedy gridding of a permutation, showing its peg
points as hollow circles; the peg permutation for this permutation
is $5431276$ while its non-peg vector is $(0,5,0,2)$.  Note that
since this is a greedy gridding, the $(1,3)$ entry of the
corresponding matrix must be $1$.} \label{fig-pegs-example}
\end{figure}

We define a {\it peg point\/} of $\pi$ to be a point which is
either first or last (either horizontally or vertically; since the
blocks are monotone, it doesn't matter) in its block in the greedy
$M$-gridding of $\pi$.  An example is shown in
Figure~\ref{fig-pegs-example}.  The {\it peg permutation\/},
$\rho^\pi$, of $\pi$ is then the permutation formed by its peg
points.  We also associate to each permutation $\pi\in\C$ its {\it
non-peg vector\/} $\y^\pi=(y_1,y_2,\dots,y_t)$, where $y_i$
denotes the number of non-peg points in
$\pi([c_i,c_{i+1})\times[n])$.  Because the $M$-gridding was
chosen greedily, the pair $(\rho^\pi, \y^\pi)$ uniquely determines
$\pi$.

We now partition the class $\C$ based upon peg permutations. Since
there can be at most $3^t$ different peg permutations of members
of $\C$ (for every column of $M$ a peg permutation can have $0$,
$1$, or $2$ elements), this is a partition into a finite number of
subsets.  Let $\C^\rho$ denote the subset of $\C$ with peg
permutation $\rho$.  This is not a permutation class (the peg
permutation of $\sigma\le\pi$ need not be the peg permutation of
$\pi$), but the set of non-peg vectors of permutations in this
class, $\{\y^\pi : \pi\in\C^\rho\}$, is a downset of vectors
in $\mathbb{N}^t$.  Therefore Theorem~\ref{vector-poly} shows that
$\C^\rho$ has eventually polynomial enumeration, and so $\C$ does
as well.
\end{proof}

\section{The Fibonacci dichotomy}\label{polynomials}

The Fibonacci dichotomy for permutation classes, first proved by
Kaiser and Klazar~\cite{kk:growth}, states that all sub-Fibonacci
permutation classes\footnote{We call a class $\C$ {\it sub-Fibonacci\/}
if $|\C_n|$ is strictly less than the $n$th Fibonacci number for some $n$.
The definition of sub-$2^{n-1}$ is analogous.}
have eventually
polynomial enumeration.
Here we give a new proof using the
characterisation of grid classes.
We have already shown, in
Theorem~\ref{grid-poly}, that grid classes of matchings and their
subclasses have eventually polynomial enumeration.  It remains
only to show that all sub-Fibonacci classes lie in grid classes of
matchings. We do this in two parts.  First we observe in
Proposition~\ref{sub-fib-griddable} that all sub-Fibonacci classes
are griddable, and then we show in
Proposition~\ref{interleavings-matchings} that all sub-$2^{n-1}$
griddable classes (which includes sub-Fibonacci classes) lie in
grid classes of matchings.

\begin{proposition}\label{sub-fib-griddable}
All sub-Fibonacci classes are griddable.
\end{proposition}
\begin{proof}
Let $\C$ denote a non-griddable class, so, by
Theorem~\ref{griddable-characterisation} and symmetry we may
assume that $\C$ contains arbitrarily long direct sums of $21$.
Since $\C$ is a permutation class, it must also contain every
permutation that embeds into an arbitarily long direct sum of
$21$.  These permutations have the form
$\sigma_1\oplus\sigma_2\oplus\cdots\oplus\sigma_k$ where each
$\sigma_i$ is either $1$ or $21$.  Thus there are precisely as
many permutations of this form of length $n$ as there are ways of
writing $n$ as an ordered sum of $1$'s and $2$'s, of which there
are $F_n$.
\end{proof}

\begin{figure}
\begin{center}
\begin{tabular}{ccc}
\psset{xunit=0.01in, yunit=0.01in} \psset{linewidth=1\psxunit}
\begin{pspicture}(0,0)(150,150)
\psline[linecolor=gray](0,75)(150,75) \psaxes[dy=10, Dy=1, dx=10,
Dx=1, tickstyle=bottom, showorigin=false,
labels=none](0,0)(140,140) \pscircle*(10,70){4.0\psxunit}
\pscircle*(20,100){4.0\psxunit} \pscircle*(30,30){4.0\psxunit}
\pscircle*(40,120){4.0\psxunit} \pscircle*(50,50){4.0\psxunit}
\pscircle*(60,80){4.0\psxunit} \pscircle*(70,60){4.0\psxunit}
\pscircle*(80,140){4.0\psxunit} \pscircle*(90,10){4.0\psxunit}
\pscircle*(100,110){4.0\psxunit} \pscircle*(110,40){4.0\psxunit}
\pscircle*(120,130){4.0\psxunit} \pscircle*(130,20){4.0\psxunit}
\pscircle*(140,90){4.0\psxunit}
\end{pspicture}
& \rule{20pt}{0pt} & \psset{xunit=0.01in, yunit=0.01in}
\psset{linewidth=1\psxunit}
\begin{pspicture}(0,0)(150,150)
\psline[linecolor=gray](75,0)(75,150) \psaxes[dy=10, Dy=1, dx=10,
Dx=1, tickstyle=bottom, showorigin=false,
labels=none](0,0)(140,140) \pscircle*(10,90){4.0\psxunit}
\pscircle*(20,130){4.0\psxunit} \pscircle*(30,30){4.0\psxunit}
\pscircle*(40,110){4.0\psxunit} \pscircle*(50,50){4.0\psxunit}
\pscircle*(60,70){4.0\psxunit} \pscircle*(70,10){4.0\psxunit}
\pscircle*(80,60){4.0\psxunit} \pscircle*(90,140){4.0\psxunit}
\pscircle*(100,20){4.0\psxunit} \pscircle*(110,100){4.0\psxunit}
\pscircle*(120,40){4.0\psxunit} \pscircle*(130,120){4.0\psxunit}
\pscircle*(140,80){4.0\psxunit}
\end{pspicture}
\end{tabular}
\end{center}
\caption{A vertical alternation (left) and its inverse, a
horizontal alternation (right).}\label{fig-alternations}
\end{figure}
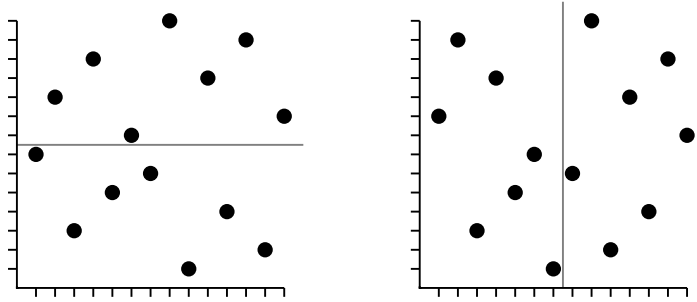

We define a {\it vertical alternation\/} to be a permutation
$\pi\in S_{2n}$ so that either
$$
\pi(1), \pi(3), \dots, \pi(2n-1)<\pi(2), \pi(4),\dots, \pi(2n)
$$
or
$$
\pi(1), \pi(3), \dots, \pi(2n-1)>\pi(2), \pi(4),\dots, \pi(2n).
$$
A {\it horizontal alternation\/} is then the group-theoretic
inverse of a vertical alternation.  Examples are shown in
Figure~\ref{fig-alternations}.  We begin by observing that classes
with arbitrarily long alternations are not small.

\begin{proposition}\label{alternations-not-small}
If the permutation class $\C$ contains arbitrarily long
alternations, then $|\C_n|\ge 2^{n-1}$ for all $n$.
\end{proposition}
\begin{proof}
By symmetry, let us suppose that $\C$ contains arbitrarily long
horizontal alternations.  By the Erd\H{o}s-Szekeres Theorem, $\C$
contains arbitrarily long horizontal alternations in which both
sides are monotone.  Therefore $\C$ contains either $\Grid(1\ 1)$,
$\Grid(1\ -1)$, $\Grid(-1\ 1)$, or $\Grid(-1\ -1)$.  It is easy
to compute that the first and last of these classes contain
$2^n-n$ permutations of length $n$ for $n\ge 1$ while the second
and third contain $2^{n-1}$ permutations of length $n\ge 1$,
establishing the proposition.
\end{proof}

Therefore a sub-Fibonacci class cannot contain arbitrarily long
alternations.  We now prove that this implies that these
classes lie in grid classes of matchings.

We say that a set of indices $\{i_1,i_2,\dots,i_s\}$ in $\pi$ is
an {\it uninterrupted monotone interval\/} if $|i_{j+1}-i_j|=1$
and $|\pi(i_{j+1})-\pi(i_j)|=1$ for all $j\in[s-1]$.  Note that if
$G(M)$ is a matching, then an $M$-gridding of $\pi$ is a division
of the elements of $\pi$ into uninterrupted monotone intervals.
Conversely, every division of $\pi$ into uninterrupted monotone
intervals gives an $M$-gridding of $\pi$ for some $M$ where $G(M)$
is a matching.

\begin{proposition}\label{interleavings-matchings}
A griddable class lies in the grid class of a matching if and only
if it does not contain arbitrarily long alternations.
\end{proposition}
\begin{proof}
One direction is obvious: if a permutation class contains
arbitrarily long alternations then it cannot lie in the grid class
of a matching.   The other direction is almost as clear, but a
formal proof takes a modest amount of effort.

Let $\C\subseteq\Grid(N)$ for some $t\times u$ matrix $N$, and
suppose that $\C$ does not contain any alternations (either
horizontal or vertical) with more than $d$ elements.  It suffices
to show that there is a constant $m$ such that every permutation
$\pi\in\C$ lies in $\Grid(M)$ where $G(M)$ is a matching and $M$
(which we allow to depend on $\pi$) has at most $m$ nonzero
entries.  This is because we can ignore the all-$0$ rows and columns,
so the size of $M$ can be bounded, and then there are only
finitely many such matrices, so $\C$ will lie in the grid class of
their direct sum (which also has a matching for its graph).
Equivalently, by our remarks above, it suffices to show that every
permutation in $\C$ can be divided into a bounded number of
uninterrupted monotone intervals.

To this end, take some permutation $\pi\in\C$ of length $n$ with
$N$-gridding given by $1=c_1\le c_2\le\cdots\le c_{t+1}=n+1$ and
$1=r_1\le r_2\le\cdots\le r_{u+1}=n+1$ and consider a particular
block in this gridding, say
$$
\pi^{(k,\ell)}:= \pi([c_k,c_{k+1})\times[r_\ell,r_{\ell+1})).
$$
We consider four types of alternations that elements of this
block can participate in: vertical alternations either with blocks
of the form $\pi^{(k,\ell^+)}$ for $\ell^+>\ell$ or of the form
$\pi^{(k,\ell^-)}$ for $\ell^-<\ell$, and horizontal alternations
with blocks of the form $\pi^{(k^+,\ell)}$ for $k^+>k$ or of the
form $\pi^{(k^-,\ell)}$ for $k^-<k$.

Every time that two consecutive elements in a block are separated
either horizontally or vertically (that is, every time that two
consecutive elements in a block fail to lie in an uninterrupted
monotone interval together) then they contribute to the length of
at least one of these four alternations.  Therefore, at most $4d$
such separations can occur, so $\pi^{(k,\ell)}$ can be divided
into at most $4d+1$ uninterrupted monotone intervals.  Hence $\pi$
itself can be divided into at most $(4d+1)tu$ uninterrupted
monotone intervals, proving the proposition.
\end{proof}

Having established that sub-$2^{n-1}$ griddable classes (and in
particular, sub-Fibonacci classes) lie in grid classes of
matchings, we now have another proof of the Fibonacci dichotomy:

\begin{corollary}\label{dichotomy}
For every permutation class $\C$, one of the following occurs:
\begin{itemize}
\item $|\C_n|\ge F_n$ for all $n$, or \item $\C$ lies in the grid
class of a matching and there is a polynomial $p(n)$ so that
$|\C_n|=p(n)$ for all sufficiently large $n$.
\end{itemize}
\end{corollary}

\section{Concluding remarks}

\noindent{\bf Decidability.}  It is not hard to see that the hypotheses of our characterisation theorems are decidable from the basis of a finitely based class.  For example, in order to determine if $\Av(B)$ contains arbitrarily long direct sums of $12$ one needs only check if any element of $B$ lies in $\Av(231,312,321)$, which is the set of permutations that are contained in arbitrarily long direct sums of $12$.  Thus we have the following result.

\begin{corollary}
Given a finite set of permutations $B$, it is decidable whether or not $\Av(B)$ is griddable.
\end{corollary}

Similar arguments show that polynomial enumeration is decidable for finitely based classes.  One first needs to check whether the class is griddable and then decide whether the class contains arbitrarily long alternations.

\begin{corollary}
Given a finite set of permutations $B$, it is decidable whether or not $|\Av_n(B)|$ agrees with a polynomial for all sufficiently large $n$.
\end{corollary}

\bigskip\noindent{\bf Finite bases.}  The decidability results above only apply to finitely based classes, however, it happens that permutation classes with polynomial enumeration must be finitely based.  Because these classes lie in grid classes of matchings, they also lie in grid classes of $\zpm$ row vectors.  Now one needs only to apply the result of Atkinson, Murphy, and Ru\v{s}kuc~\cite{amr:pwocsop} that every subclass of the grid class of a $\zpm$ row vector is finitely based.

\bigskip\noindent{\bf Enumeration.}  While Corollary~\ref{dichotomy} characterises the sub-Fibonacci
classes and shows that they have eventually polynomial
enumeration, it does not address the issue of how one might find
these formulas.  This could presumably be settled
by strengthening the results given here to obtain bounds
(computable from the basis of $\C$) on the degree of the
polynomial and the values of $n$ for which $|\C_n|$ agrees with
this polynomial, but there are already three general methods which can be used to count these classes:
\begin{enumerate}
\item[(i)] Since permutation classes with polynomial growth lie in grid classes of matchings, they also lie in grid classes of $\zpm$ row vectors.  One can therefore use Theorem~\ref{star-rational} to enumerate them.
\item[(ii)] Permutation classes with polynomial growth contain only finitely many simple permutations, and so the results of Albert and Atkinson~\cite{aa:simple:alg} apply to them.
\item[(iii)] Permutation classes with polynomial growth correspond to regular languages via the insertion encoding of Albert, Linton, and Ru\v{s}kuc~\cite{insertion}, so this method can also be used.
\end{enumerate}

\bigskip\noindent{\bf Classes with two basis elements.}  Albert, Atkinson, and Brignall~\cite{aab:poly} have recently studied doubleton-based classes with polynomial enumeration.  They give a characterisation of the pairs of permutations $\beta_1$ and $\beta_2$ for which $\Av(\beta_1,\beta_2)$ has polynomial enumeration and give bounds on the degree of this polynomial.

\bigskip

\noindent{\it Acknowledgment.}  We thank Nik Ru\v{s}kuc and Bruce Sagan for their helpful comments.

\bibliographystyle{acm}
\begin{small}
\bibliography{refs}

\def\cprime{$'$}
\begin{thebibliography}{10}

\bibitem{aa:simple:alg}
{\sc Albert, M.~H., and Atkinson, M.~D.}
\newblock {Simple permutations and pattern restricted permutations}.
\newblock {\em Discrete Math. 300}, 1-3 (2005), 1--15.

\bibitem{aab:poly}
{\sc Albert, M.~H., Atkinson, M.~D., and Brignall, R.}
\newblock Permutation classes of polynomial growth.
\newblock In preparation.

\bibitem{aar:regular}
{\sc Albert, M.~H., Atkinson, M.~D., and Ru{\v{s}}kuc, N.}
\newblock Regular closed sets of permutations.
\newblock {\em Theoret. Comput. Sci. 306}, 1-3 (2003), 85--100.

\bibitem{insertion}
{\sc Albert, M.~H., Linton, S., and Ru\v{s}kuc, N.}
\newblock The insertion encoding of permutations.
\newblock {\em Electron. J. Combin. 12}, 1 (2005), Research paper 47, 31 pp.
  (electronic).

\bibitem{a:skewmerged}
{\sc Atkinson, M.~D.}
\newblock Permutations which are the union of an increasing and a decreasing
  subsequence.
\newblock {\em Electron. J. Combin. 5\/} (1998), Research paper 6, 13 pp.\
  (electronic).

\bibitem{a:rp}
{\sc Atkinson, M.~D.}
\newblock Restricted permutations.
\newblock {\em Discrete Math. 195}, 1-3 (1999), 27--38.

\bibitem{amr:pwocsop}
{\sc Atkinson, M.~D., Murphy, M.~M., and Ru{\v{s}}kuc, N.}
\newblock Partially well-ordered closed sets of permutations.
\newblock {\em Order 19}, 2 (2002), 101--113.

\bibitem{bbw:growth2000}
{\sc Balogh, J., Bollob{\'a}s, B., and Weinreich, D.}
\newblock The speed of hereditary properties of graphs.
\newblock {\em J. Combin. Theory Ser. B 79}, 2 (2000), 131--156.

\bibitem{es:acpig}
{\sc Erd\H{o}s, P., and Szekeres, G.}
\newblock A combinatorial problem in geometry.
\newblock {\em Compos. Math. 2\/} (1935), 463--470.

\bibitem{fh:sg}
{\sc F{\"o}ldes, S., and Hammer, P.~L.}
\newblock Split graphs.
\newblock In {\em Proceedings of the Eighth Southeastern Conference on
  Combinatorics, Graph Theory and Computing (Louisiana State Univ., Baton
  Rouge, La., 1977)\/} (Winnipeg, Man., 1977), Utilitas Math., pp.~311--315.
  Congressus Numerantium, No. XIX.

\bibitem{kk:growth}
{\sc Kaiser, T., and Klazar, M.}
\newblock On growth rates of closed permutation classes.
\newblock {\em Electron. J. Combin. 9}, 2 (2002/03), Research paper 10, 20 pp.
  (electronic).

\bibitem{ksw:incdec}
{\sc K{\'e}zdy, A.~E., Snevily, H.~S., and Wang, C.}
\newblock Partitioning permutations into increasing and decreasing
  subsequences.
\newblock {\em J. Combin. Theory Ser. A 73}, 2 (1996), 353--359.

\bibitem{knuth1}
{\sc Knuth, D.~E.}
\newblock {\em The art of computer programming. {V}ol. 1: {F}undamental
  algorithms}.
\newblock Addison-Wesley Publishing Co., Reading, Mass., 1969.

\bibitem{ls:smooth}
{\sc Lakshmibai, V., and Sandhya, B.}
\newblock Criterion for smoothness of {S}chubert varieties in {${\rm
  Sl}(n)/B$}.
\newblock {\em Proc. Indian Acad. Sci. Math. Sci. 100}, 1 (1990), 45--52.

\bibitem{mt:swc}
{\sc Marcus, A., and Tardos, G.}
\newblock Excluded permutation matrices and the {S}tanley-{W}ilf conjecture.
\newblock {\em J. Combin. Theory Ser. A 107}, 1 (2004), 153--160.

\bibitem{profile}
{\sc Murphy, M.~M., and Vatter, V.}
\newblock Profile classes and partial well-order for permutations.
\newblock {\em Electron. J. Combin. 9}, 2 (2002/03), Research paper 17, 30 pp.
  (electronic).

\bibitem{pt:poly}
{\sc Pouzet, M., and Thi{\'e}ry, N.~M.}
\newblock Some relational structures with polynomial growth and their
  associated algebras.
\newblock arXiv:math.CO/0601256.

\bibitem{sz:hered}
{\sc Scheinerman, E.~R., and Zito, J.}
\newblock On the size of hereditary classes of graphs.
\newblock {\em J. Combin. Theory Ser. B 61}, 1 (1994), 16--39.

\bibitem{stankova:fs}
{\sc Stankova, Z.~E.}
\newblock Forbidden subsequences.
\newblock {\em Discrete Math. 132}, 1-3 (1994), 291--316.

\bibitem{s:E2546}
{\sc Stanley, R.~P.}
\newblock Solution to problem {E}2546.
\newblock {\em Amer. Math. Monthly 83}, 10 (1976), 813--814.

\end{thebibliography}
\end{small}

\end{document}